\documentclass{ifacconf}
\usepackage{amsmath,amsfonts,amssymb}
\usepackage{dsfont}
\usepackage{graphicx}      
\usepackage{natbib}
\usepackage{bbm}    
\usepackage{mathrsfs}
\newtheorem{theorem}{Theorem}[section]
\newtheorem{lemma}[theorem]{Lemma}   
\newtheorem{remark}[theorem]{Remark}   
\newtheorem{proposition}[theorem]{Proposition} 
\newtheorem{conjecture}[theorem]{Conjecture}

\begin{document}
\begin{frontmatter}

\title{On estimation and feedback control of spin-$\frac12$ systems with unknown initial states\thanksref{footnoteinfo}} 

\thanks[footnoteinfo]{The authors thank the support from Agence Nationale de la Recherche projects Q-COAST ANR-19-CE48-0003 and QUACO ANR-17-CE40-0007.}

\author[First]{Weichao Liang} 
\author[First]{Nina H. Amini} 
\author[First]{Paolo Mason}

\address[First]{Laboratoire des Signaux et Syst\` emes, CNRS - CentraleSup\'elec - Univ. Paris-Sud, Universit\'e Paris-Saclay, 3, rue Joliot Curie, 91192, Gif-sur-Yvette, France.  (e-mail: first name.family name@centralesupelec.fr).}

\begin{abstract}                
In this paper, we consider stochastic master equations describing the evolutions of quantum systems interacting with electromagnetic fields undergoing continuous-time measurements. In particular, we study feedback control of quantum spin-$\frac12$ systems in the case of unawareness of initial states and in presence of measurement imperfections. We prove that the fidelity between the actual quantum filter  and  its associated estimated filter converges to one under appropriate assumption on the feedback controller. This shows the asymptotic convergence of such filters. In addition, for spin-$J$ systems, we discuss heuristically the asymptotic behavior of the actual quantum filter and its associated estimated filter and the possibility of  exponentially stabilizing such systems towards an eigenvector of the measurement  operator by an appropriate feedback. 
\end{abstract}

\begin{keyword}
Spin systems; Quantum filtering; Stochastic master equations; Fidelity.
\end{keyword}
\end{frontmatter}
\section{Introduction}
Classical filtering~\cite{kallianpur2013stochastic,xiong2008introduction} determines
the best estimation of the state of a classical system from noisy observations, the quantum
analogue was developed in the 1960s by~\cite{davies1969quantum,davies1976quantum} and extended by Belavkin
in the 1980s~\cite{belavkin1983theory,belavkin1989nondemolition,belavkin1992quantum,belavkin1995quantum}, relying on the quantum probability theory
and the quantum stochastic calculus~\cite{hudson1984quantum,hudson2003introduction,meyer2006quantum}. The modern treatment of
quantum filtering has been established in~\cite{bouten2009discrete}.
Roughly speaking, quantum filtering theory gives  a matrix-valued
stochastic differential equation called stochastic master equation, to describe the time
evolution of the state of an open quantum system interacting with an electromagnetic field under homodyne detection.

In real experiments, different types of imperfections, such as detection inefficiencies and unawareness of initial states, may be present (see e.g.,~\cite{sayrin2011real}).  In the case of unawareness of initial states, one considers estimated quantum filters which are designed based on measurements. The main question is whether the filter ``forgets" its initial state and has the same asymptotic as the actual quantum filter. This problem can be posed when a feedback depending on the estimated quantum filter is applied. Since the observation process depends on the actual filter state, we deal with coupled stochastic master equations whose asymptotic behavior is not at all a trivial problem.


When the feedback is turned off, this problem has been investigated in some recent papers. In the series of papers~\cite{van2006filtering,van2009observability,van2009stability,van2010nonlinear}, a sufficient observability condition has been established so that such convergence is guaranteed. However, such condition is not easy to verify even if the system is finite dimensional. In~\cite{diosi2006coupled}, the problem of convergence is addressed assuming  either the actual quantum filter or its estimated state is always pure. The authors prove that the fidelity is  a sub-martingale for this case. Then, in~\cite{rouchon2011fidelity}, by applying Uhlmann's technique~\cite[Theorem 9.4]{nielsen2002quantum}, the author shows that the fidelity between the state of the discrete-time quantum filter and its associated estimated state is a sub-martingale via a Kraus map. However, this sub-martingale property of the fidelity cannot ensure the convergence of the filter state towards the actual one. In~\cite{amini2011stability}, the authors show that the fidelity is a sub-martingale for continuous-time quantum filters with perfect measurement for arbitrary mixed states. By the quantum repeated interaction approach, such result has been extended to the continuous-time jump-diffusion stochastic master equations  with general measurement imperfections in~\cite{amini2014stability}. Then, in~\cite{benoist2014large}, for quantum filters described by jump-diffusion stochastic differential equations,  the authors show that when the control input is turned off, under perfect Quantum Non-Demolition (QND) measurements and a non-degeneracy assumption, the convergence is ensured.

Concerning the feedback stabilization of discrete-time QND measures, in~\cite{mirrahimi2009feedback}, the authors show the convergence towards a Fock state in the case of unawareness of  initial states. This has been obtained under appropriate assumptions on the initial states of the filter and its estimate. For continuous-time case, stabilization results for the case of unknown initial states have not been addressed so far. On the other hand, for angular momentum systems and without initialization imperfections, in~\cite{mirrahimi2007stabilizing}, the authors show the asymptotic convergence towards a chosen eigenvector of the measurement operator $J_z$. In the same context, in~\cite{liang2019exponential}, we provide general conditions on the feedback controller and a local Lyapunov type condition which ensure exponential convergence towards a chosen eigenvector of $J_z.$ 

In this paper, we first analyze the dynamics of quantum spin-$\frac12$ systems in presence of feedback control. We suppose imperfections in measurements and unawareness of the initial state. We show that the filter and filter estimate have the same asymptotic behavior under appropriate assumption on the feedback. For spin-$J$ systems, we discuss heuristically the asymptotic behavior of the actual quantum filter and its associated estimated filter and the possibility of  exponentially stabilizing such systems towards an eigenvector of the measurement  operator $J_z$  by a candidate feedback controller. Numerical simulations are provided in order to illustrate our results and to support the efficiency of the proposed candidate feedback.
\section{Model description}
Here, we consider quantum spin-$\frac{1}{2}$ systems. The stochastic master equations describing the evolution of the actual system state and the corresponding estimated state are given as follows,
\begin{align}
d\rho_t&\!=\!F_{\hat u_t}(\rho_t)dt\!+\!L(\rho_t)dt\!+\!G(\rho_t)\big(dY_t\!-\!2\sqrt{\eta M}\mathrm{Tr}(\sigma_z \rho_t)dt\big),\nonumber\\
d\hat{\rho}_t&\!=\!F_{\hat u_t}(\hat{\rho}_t)dt\!+\!L(\hat{\rho}_t)dt\!+\!G(\hat{\rho}_t)\big(dY_t\!-\!2\sqrt{\eta M}\mathrm{Tr}(\sigma_z \hat{\rho}_t)dt\big),\nonumber
\end{align}
where
\begin{itemize}
\item the actual quantum state of the spin-$\frac12$ system is denoted as $\rho$, and belongs to the space $
\mathcal{S}_2:=\{\rho\in\mathbb{C}^{2\times 2}|\,\rho=\rho^*,\mathrm{Tr}(\rho)=1,\rho \geq 0\}$. The associated estimated state is denoted as $\hat{\rho}\in\mathcal{S}_2$,
\item the matrices $\sigma_x,$ $\sigma_y$ and $\sigma_z$ correspond to the Pauli matrices. 
\item $F_{u}(\rho):=-i[\omega\sigma_z+u\sigma_y,\rho]$, $L(\rho):=M(\sigma_z\rho\sigma_z-\rho)$ and $G(\rho):=\sqrt{\eta M}\big(\sigma_z\rho+\rho\sigma_z-2\mathrm{Tr}(\sigma_z\rho)\rho\big)$. 
\item $Y_t$ denotes the observation process of the actual quantum spin-$\frac12$ system, which is a continuous semi-martingale whose quadratic variation is given by $[Y,Y]_t=t$. Its dynamics satisfies $dY_t=dW_t+2\sqrt{\eta M}\mathrm{Tr}(\sigma_z\rho_t)dt$, where $W_t$ is a one-dimensional standard Wiener process,
\item $\hat u_t:=u(\hat{\rho}_t)$ denotes the feedback controller as a function of the estimated state $\hat{\rho}_t$,
\item $\omega$ is the difference between the energies of the excited state and the ground state, $\eta\in[0,1]$ describes the efficiency of the detector, and $M>0$ is the strength of the interaction between the system and the probe. 
\end{itemize}
By replacing $dY_t=dW_t+2\sqrt{\eta M}\mathrm{Tr}(\sigma_z\rho_t)dt$ in the equation above, we obtain the following matrix-valued stochastic differential equations describing the time evolution of the pair $(\rho_t,\hat{\rho}_t)\in\mathcal{S}_2\times\mathcal{S}_2$,
 \begin{align}
 d\rho_t&=F_{\hat u_t}(\rho_t)dt+L(\rho_t)dt+G(\rho_t)dW_t,\label{Ch5_StabFilter:2D SME}\\
 d\hat{\rho}_t&=F_{\hat u_t}(\hat{\rho}_t)dt+L(\hat{\rho}_t)dt+2\sqrt{\eta M}G(\hat{\rho}_t)\mathrm{Tr}\big(\sigma_z(\rho_t-\hat{\rho}_t)\big)dt\nonumber\\
 &\quad+G(\hat{\rho}_t)dW_t.\label{Ch5_StabFilter:2D SME filter}
 \end{align}
If $u\in\mathcal{C}^1(\mathcal{S}_2,\mathbb{R})$, the existence and uniqueness of the solution of~\eqref{Ch5_StabFilter:2D SME}--\eqref{Ch5_StabFilter:2D SME filter} can be shown by similar arguments as in~\cite[Proposition 3.5]{mirrahimi2007stabilizing}.
Recall that a density operator can be uniquely characterized by the Bloch sphere coordinates $(x,y,z)$ as 
\begin{equation*}
\rho=\frac{\mathds{1}+x\sigma_x+y\sigma_y+z\sigma_z}{2}=\frac 12
\begin{bmatrix}
1+z & x-iy\\
x+iy & 1-z
\end{bmatrix}.
\end{equation*}
The vector $(x,y,z)$ belongs to the ball
\begin{equation*}
\mathcal{B}:=\{(x,y,z)\in\mathbb{R}^3|\, x^2+y^2+z^2 \leq 1 \}.
\end{equation*}
The stochastic differential equation~\eqref{Ch5_StabFilter:2D SME} expressed in the Bloch sphere coordinates takes the following form
\begin{subequations}
\begin{align}
dx_t&=\left(\!-\omega_{eg} y_t-\frac{M}{2}x_t+\hat{u}_tz_t \!\right)dt\!-\!\sqrt{\eta M}x_tz_tdW_t,\label{SME Bloch:x}\\
dy_t&=\left(\omega_{eg} x_t-\frac{M}{2} y_t \right)dt-\sqrt{\eta M}y_tz_tdW_t, \label{SME Bloch:y}\\
dz_t&=-\hat{u}_tx_tdt+\sqrt{\eta M} ( 1-z^2_t )dW_t.\label{SME Bloch:z}
\end{align}
\label{SME Bloch}
\end{subequations}
The stochastic differential equation~\eqref{Ch5_StabFilter:2D SME filter} in the Bloch sphere coordinates is given by,
\begin{subequations}
\begin{align}
d\hat{x}_t=&\left(-\omega_{eg} \hat{y}_t-\frac{M}{2}\hat{x}_t+\hat{u}_t\hat{z}_t+\eta M\hat{x}_t\hat{z}_t(\hat{z}_t-z_t) \right)dt\nonumber\\
&-\sqrt{\eta M}\hat{x}_t\hat{z}_tdW_t,\label{SME filter Bloch:x}\\
d\hat{y}_t=&\left(\omega_{eg} \hat{x}_t-\frac{M}{2} \hat{y}_t +\eta M\hat{y}_t\hat{z}_t(\hat{z}_t-z_t) \right)dt\nonumber\\
&-\sqrt{\eta M}\hat{y}_t\hat{z}_tdW_t, \label{SME filter Bloch:y}\\
d\hat{z}_t=&\left(-\hat{u}_t\hat{x}_t-\eta M(1-\hat{z}^2_t)(\hat{z}_t-z_t) \right)dt\nonumber\\
&+\sqrt{\eta M}( 1-\hat{z}^2_t)dW_t.\label{SME filter Bloch:z}
\end{align}
\label{SME filter Bloch}
\end{subequations}
\section{Convergence property of quantum spin-$\frac12$ systems}
We focus on the fidelity $\mathcal{F}(\rho,\hat{\rho})$ which defines a ``distance" between the real state $\rho$ and the estimated state $\hat{\rho}$. In the two-level case, the fidelity can be written in the following form
\begin{equation*}
\mathcal{F}(\hat{\rho},\rho)=\mathrm{Tr}(\hat{\rho}\rho)+2\sqrt{\det(\hat{\rho})\det(\rho)}.
\end{equation*}
Thus the fidelity in the Bloch sphere coordinates is given by
\begin{equation*}
\mathcal{F}(\hat{\rho},\rho)\!=\!\mathcal{F}(\hat{\mathbf{v}},\mathbf{v})\!=\!\frac{1}{2}\!\left(1\!+\!\mathbf{v}^{\top}\hat{\mathbf{v}}\!+\!\sqrt{(1-\|\mathbf{v}\|^2)(1-\|\hat{\mathbf{v}}\|^2)}  \right)\!,
\end{equation*}
where $\mathbf{v}:=(x,y,z)$ denotes the real state and $\hat{\mathbf{v}}:=(\hat{x},\hat{y},\hat{z})$ denotes the estimated state in Bloch sphere coordinates. Thus, for the two special cases $\mathcal{F}(\hat{\rho},\rho)\!=\!1$ and $\mathcal{F}(\hat{\rho},\rho)\!=\!0$,
\begin{enumerate}
\item if $\mathcal{F}(\hat{\rho},\rho)=1$, we have $\mathbf{v}=\hat{\mathbf{v}}$;
\item if $\mathcal{F}(\hat{\rho},\rho)=0$, we have $\mathbf{v}+\hat{\mathbf{v}}=0$ and $\|\mathbf{v}\|^2\!=\!\|\hat{\mathbf{v}}\|^2\!=\!1$.
\end{enumerate}
In order to apply the It\^o formula on the fidelity $\mathcal{F}(\rho,\hat{\rho})$, we need to show the unattainability of the boundary for $\rho$ and $\hat{\rho}$. By straightforward calculations, we can show that 
\begin{equation}
\{\rho\in\mathcal{S}_2|\,\det(\rho)=0\}=\{\rho\in\mathcal{S}_2|\,\mathrm{Tr}(\rho^2)=1\},
\label{Ch5_StabFilter:BoudaryPurestates}
\end{equation}
which means that the boundary $\partial\mathcal{S}_2$ is equal to the set of all pure states $\mathcal{P}$. 
The following lemma states some invariance properties for Equations~\eqref{Ch5_StabFilter:2D SME}--\eqref{Ch5_StabFilter:2D SME filter}. 
\begin{lemma}
If $\rho_0>0$, then $\mathbb{P}(\rho_t>0,\,\forall t\geq0)=1$. Moreover, if $\eta=1$, $\partial \mathcal{S}_2\times \mathcal{S}_2$ is a.s. invariant for Equations~\eqref{Ch5_StabFilter:2D SME}--\eqref{Ch5_StabFilter:2D SME filter}. The same results hold true for $\mathcal{S}_2\times \partial \mathcal{S}_2$.
\label{Ch5_StabFilter:Invaraint}
\end{lemma}
\begin{pf}
The dynamics of the purification function $S(\rho_t):=1-\mathrm{Tr}(\rho^2_t)$ is given by
\begin{align*}
dS(\rho_t)=&M\big(\frac{(1-\eta)(1-z^2_t)}{2}-(1-\eta z^2_t)S(\rho_t)\big)dt\\
&-2\sqrt{\eta M}z_tS(\rho_t)dW_t.
\end{align*}
Then, if $\eta=1$, it is obvious that the set of all pure states $\mathcal{P}$ for Equation~\eqref{Ch5_StabFilter:2D SME} is a.s. invariant.

Next, let us prove the first part of the lemma. Given $\varepsilon>0$, consider any $\mathcal{C}^2$ function on $\mathcal{S}$ such that
\begin{equation*}
V(\rho)=\frac{1}{S(\rho)},\quad \text{if }S(\rho)>\varepsilon.
\end{equation*}
We find
\begin{align*}
\mathscr{L}V(\rho)&=M\left( 1+3\eta z^2- (1-\eta)\frac{1-z^2}{2S(\rho)} \right)V(\rho),\\
&\leq 4MV(\rho)\quad\text{if }S(\rho)>\varepsilon.
\end{align*}
To conclude the proof, one applies standard arguments (see e.g., \cite[Lemma~4.1]{liang2019exponential}). Roughly speaking, by setting $f(\rho,t)=e^{-4M t} V(\rho)$, one has $\mathscr{L}f\leq 0$ whenever $S(\rho)>\varepsilon$. From this fact one proves that the probability of $S(\rho)$ becoming zero in a finite fixed time $T$ is proportional to $\varepsilon$ and, being the latter arbitrary, it must be $0$. Due to the equality~\eqref{Ch5_StabFilter:BoudaryPurestates}, $\mathbb{P}(\rho_t>0,\,\forall t\geq0)=1$ when $\rho_0>0$. The last part of the lemma can be proved in the same manner.
\end{pf}
\begin{remark}
Note that the variation of the purification function does not explicitly depend on the feedback.  
\end{remark}
Next, we analyze the behavior of $\rho_t$ (resp. $\hat{\rho}_t$) whenever the corresponding initial datum $\rho_0$ (resp. $\hat{\rho}_0$) lies at the boundary of $S_2$.
Denote $\boldsymbol \rho_g:=\mathrm{diag}(1,0)$ and $\boldsymbol \rho_e:=\mathrm{diag}(0,1)$, which are the the pure states corresponding to the eigenvectors of~$\sigma_z$. 
\begin{lemma}
Assume $\eta\in(0,1)$ and $u\in\mathcal{C}^1(\mathcal{S}_2,\mathbb{R})$. Suppose that $\hat{\rho}_0$ lies in $\partial\mathcal{S}_2\setminus\{\boldsymbol \rho_g, \boldsymbol \rho_e\},$ then $\hat{\rho}_t>0$ for all $t>0$ almost surely. Moreover, if $\hat{\rho}_0\in\{\boldsymbol \rho_g, \boldsymbol \rho_e\}$ and $u(\hat\rho_0)\neq 0$ then, $\hat{\rho}_t>0$ for all $t>0$ almost surely. In addition, under the assumption $\rho_0\in\partial\mathcal{S}_2\setminus\{\boldsymbol \rho_g, \boldsymbol \rho_e\}$, then $\rho_t>0$ for all $t>0$ almost surely. Also, if $u(\boldsymbol \rho_g)u(\boldsymbol \rho_e)\neq 0$ then, for all $\rho_0\in\partial\mathcal{S}_2$, $\rho_t$ exits the boundary in finite time and stays in the interior of $\mathcal{S}_2$ almost surely.
\label{Ch5_StabFilter:PassLemmaSpin1/2}
\end{lemma}
\begin{pf}
First, consider the purification function $S(\hat\rho):=1-\mathrm{Tr}(\hat\rho^2)$ for Equation~\eqref{Ch5_StabFilter:2D SME filter}, whose dynamics is given by
\begin{align*}
dS(\hat{\rho}_t)=&M\Big(\frac{(1-\eta)(1-\hat{z}^2_t)}{2}-(1-\eta \hat{z}^2_t)S(\hat\rho_t)\\
&\,-4\eta(z_t-\hat{z}_t)\hat{z}_tS(\hat\rho_t)\Big)dt-2\sqrt{\eta M}\hat{z}_tS(\hat\rho_t)dW_t.
\end{align*}
Now, assume $\hat\rho_0\in \partial\mathcal{S}_2\setminus\{B_{\epsilon}(\boldsymbol \rho_g)\cup B_{\epsilon}(\boldsymbol \rho_e)\}\}$. By compactness, there exists a $\zeta>0$ such that $\frac{1}{2}M(1-\eta)(1-\hat{z}^2)\geq \zeta.$ Define $\tau:=\inf\{t> 0|\,\hat{\rho}_t\notin \partial\mathcal{S}_2\setminus\{B_{\epsilon}(\boldsymbol \rho_g)\cup B_{\epsilon}(\boldsymbol \rho_e)\}$, for all $\hat{\rho}_0\in \partial\mathcal{S}_2\setminus\{B_{\epsilon}(\boldsymbol \rho_g)\cup B_{\epsilon}(\boldsymbol \rho_e)\}$ and $t>0$, by Ito formula, we have
\begin{equation*}
\mathbb E(S(\hat{\rho}_{t\wedge \tau}))=\mathbb E\left(\int^{t\wedge \tau}_0\frac{1}{2}M(1-\eta)(1-\hat{z}^2_s)ds\right)\geq \zeta\mathbb E(t\wedge \tau).
\end{equation*}
By continuity and the definition of $\tau$, $S(\hat{\rho}_{t\wedge \tau})=0$ almost surely. This implies that $\mathbb E(t\wedge \tau)=0.$ Since we have $\mathbb E(t\wedge \tau)\geq t\mathbb P(\tau\geq t)$ we deduce that $\mathbb P(\tau\geq t)=0$ for all $t>0.$
Due to the arbitrariness of $\epsilon$, if $\hat\rho_0\in\partial\mathcal{S}_2\setminus\{\boldsymbol \rho_g, \boldsymbol \rho_e\}$ then $\hat{\rho}_t$ exits the boundary immediately. Combining with the strong Markov property and Lemma~\ref{Ch5_StabFilter:Invaraint}, $\hat{\rho}_t>0$ for all $t>0$, almost surely. Moreover, if $\hat\rho_0\in\{\boldsymbol \rho_g, \boldsymbol \rho_e\}$ then by the condition $u(\hat\rho_0)\neq 0$ we deduce the same result. 

For the case $\rho_0\in\partial\mathcal{S}_2\setminus\{\boldsymbol \rho_g, \boldsymbol \rho_e\}$, the above arguments can be repeated so that $\rho_t>0$ for all $t>0$ almost surely. Moreover, if $u(\boldsymbol \rho_g)u(\boldsymbol \rho_e)\neq 0$ and $u\in\mathcal{C}^1(\mathcal{S}_2,\mathbb{R})$, then there exists a neighborhood of $\boldsymbol \rho_e$ denoted by $B_{r_e}(\boldsymbol \rho_e)$ and a neighborhood of $\boldsymbol \rho_g$ denoted by $B_{r_g}(\boldsymbol \rho_g)$ such that, $u(\hat{\rho})\neq 0$ for all $\rho\in B_{r_e}(\boldsymbol \rho_e)\cup B_{r_g}(\boldsymbol \rho_g)$. By applying the same arguments  as in~\cite[Lemma~6.1]{liang2019exponential}, we can show that if $\rho_0\in\partial\mathcal{S}_2,$ $ \hat{\rho}_t$ enters in $B_{r_e}(\boldsymbol \rho_e)\cup B_{r_g}(\boldsymbol \rho_g)$ in finite time almost surely, which means that $u(\hat{\rho}_t)$ becomes non-zero in finite time almost surely. As a consequence,  $\rho_t$ exits the boundary and stays in the interior of $\mathcal{S}$ almost surely. The proof is then complete. \end{pf}
\begin{proposition}\label{prop:fidelity}
Assume $\eta\in (0,1]$ and let $u\in\mathcal{C}^1(\mathcal{S}_2,\mathbb{R})$. Then for all $(\rho_0,\hat{\rho}_0)\in \mathcal{S}_2\times\mathcal{S}_2$, either $\hat{\rho}_t$ converges to $\{\boldsymbol \rho_g, \boldsymbol \rho_e\}$ or $\mathcal{F}(\rho_t,\hat{\rho}_t)$ converges to one, almost surely. In particular, if $u(\boldsymbol \rho_e)u(\boldsymbol \rho_g)\neq 0$ then $\mathcal{F}(\rho_t,\hat{\rho}_t)$ converges to one almost surely.
\end{proposition}
\begin{pf}
In order to study the asymptotic behavior of $\mathcal{F}(\rho,\hat{\rho})$, in the following we will apply It\^o formula. For this purpose, we first need to show the $\mathcal{C}^2$ regularity of $\mathcal F$ on appropriate invariant sets.  If $\eta=1,$ then by Lemma~\ref{Ch5_StabFilter:Invaraint}, $\partial \mathcal{S}_2\times\partial \mathcal{S}_2$, $\partial \mathcal{S}_2\times\mathrm {int}(\mathcal{S}_2),$ $\mathrm {int}(\mathcal{S}_2)\times\partial \mathcal{S}_2$  and $\mathrm {int}(\mathcal{S}_2)\times\mathrm {int}(\mathcal{S}_2)$ are invariant for the coupled system~\eqref{Ch5_StabFilter:2D SME}-\eqref{Ch5_StabFilter:2D SME filter} almost surely. Moreover, if either $\rho$ or $\hat{\rho}$ belongs to the boundary of $\mathcal{S}_2$, then the fidelity takes the form $\mathcal{F}(\rho,\hat{\rho})=\mathrm{Tr}(\rho\hat{\rho}),$ which is a $\mathcal{C}^2$ function. For the case $\eta\in(0,1)$, under the assumptions that $u\in\mathcal{C}^1(\mathcal{S}_2,\mathbb{R})$ and $u(\boldsymbol \rho_e)u(\boldsymbol \rho_g)\neq 0$, by Lemma~\ref{Ch5_StabFilter:PassLemmaSpin1/2}, $(\rho_t,\hat{\rho}_t)$ exits the boundary in finite time and stays in $\mathrm {int}(\mathcal{S}_2)\times\mathrm {int}(\mathcal{S}_2)$ afterwards almost surely. Note that the fidelity function is $\mathcal C^2$ in $\mathrm {int}(\mathcal{S}_2)\times\mathrm {int}(\mathcal{S}_2).$ 

Consider the Lyapunov function $\mathcal V(\rho,\hat\rho):=1-\mathcal{F}(\rho,\hat{\rho}).$ Denote $\Xi:=\sqrt{(1-\|\mathbf{v}\|^2)(1-\|\hat{\mathbf{v}}\|^2)}$. For any $u\in\mathcal{C}^1(\mathcal{S}_2,\mathbb{R})$, the infinitesimal generator of $\mathcal{F}(\rho,\hat{\rho})$ is given by
{\small\begin{equation}
\setlength{\abovedisplayskip}{3pt}
\setlength{\belowdisplayskip}{3pt}
\begin{split}
\mathscr{L}&\mathcal{F}(\rho,\hat{\rho})= 
\frac{M(1-\eta)}{4 \Xi}\Big((1-\hat{z}^2)(1-\|\mathbf{v}\|^2)\\
&+(1-z^2)(1-\|\hat{\mathbf{v}}\|^2)
+2\hat{z}^2(1-\mathbf{v}^{\top}\hat{\mathbf{v}}-\Xi)\Xi-2(1-z\hat{z})\Xi \Big)\\
&+\frac{M}{2}(1-\hat{z}^2)(1-\mathbf{v}^{\top}\hat{\mathbf{v}}-\Xi).
\label{LFidelity}
\end{split}
\end{equation}}
In particular, if $\eta=1$, we have
{\small\begin{equation}
\setlength{\abovedisplayskip}{3pt}
\setlength{\belowdisplayskip}{3pt}
\begin{split}
\mathscr{L}\mathcal{F}(\rho,\hat{\rho})&=\frac{M}{2}(1-\hat{z}^2)(1-\mathbf{v}^{\top}\hat{\mathbf{v}}-\Xi)\\
&=M(1-\hat{z}^2)\big(1-\mathcal{F}(\rho,\hat{\rho})\big).
\label{eq:inft}
\end{split}
\end{equation}}
For $\eta=0$, we have
{\small\begin{equation}
\setlength{\abovedisplayskip}{3pt}
\setlength{\belowdisplayskip}{3pt}
\begin{split}
\mathscr{L}\mathcal{F}(\rho,\hat{\rho})=&\frac{M}{2}\bigg(\frac{(1-\hat{z}^2)(1-\|\mathbf{v}\|^2)+(1-z^2)(1-\|\hat{\mathbf{v}}\|^2)}{2 \Xi}\nonumber\\
&+z\hat{z}-\mathbf{v}^{\top}\hat{\mathbf{v}}-\Xi\bigg)\nonumber\\
\geq&\frac{M}{2}\left(\frac{(1-\hat{z}^2)(1-\|\mathbf{v}\|^2)+(1-z^2)(1-\|\hat{\mathbf{v}}\|^2)}{2 \Xi}\right.\nonumber\\
&-\sqrt{(\|\mathbf{v}\|^2-z^2)(\|\hat{\mathbf{v}}\|^2-\hat{z}^2)}-\Xi\bigg)\nonumber\\
=&\frac{M}{4\Xi}\Big(\sqrt{(\|\hat{\mathbf{v}}\|^2-\hat{z}^2)(1-\|\mathbf{v}\|^2)}\nonumber\\
&-\sqrt{(\|\mathbf{v}\|^2-z^2)(1-\|\hat{\mathbf{v}}\|^2)}\Big)^2.
\label{eq:inft0}
\end{split}
\end{equation}}
Therefore, for all $\eta\in[0,1]$ and $(\rho,\hat{\rho})\in\mathrm {int}(\mathcal{S}_2)\times\mathrm {int}(\mathcal{S}_2)$, we have $\mathscr{L}\mathcal{F}(\rho,\hat{\rho})\geq0$ which implies that $\mathscr{L}\mathcal{V}(\rho,\hat{\rho})\leq0$. By the stochastic LaSalle-type theorem in~\cite{mao1999stochastic}, we deduce that $\lim_{t\rightarrow\infty}\mathscr{L} \mathcal F(\rho_t,\hat{\rho}_t)=0$ almost surely. Since $\mathscr{L} \mathcal F(\rho,\hat{\rho})$ for any $\eta\in(0,1]$ can be written as a convex combination of the expressions~\eqref{eq:inft} and~\eqref{eq:inft0}, we have that either $|\hat z|$ converges to one or $\mathcal F(\rho,\hat{\rho})$ converges to one almost surely. 
This concludes the proof of the first part of the proposition. 
The additional assumption $u(\boldsymbol \rho_e)u(\boldsymbol \rho_g)\neq 0$ rules out the first possibility, completing the proof of the proposition.
\end{pf}
\subsection{Simulations}
\begin{figure}[thpb]
\centering
\includegraphics[width=6cm]{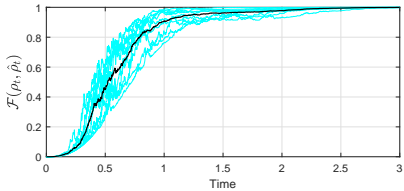}
\caption{\small Convergence of the fidelity $\mathcal{F}(\rho_t,\hat{\rho}_t)$ towards one with the feedback law $u(\hat{\rho})\equiv 1$ starting at $(\rho_0,\hat{\rho}_0)=(\boldsymbol\rho_e,\boldsymbol\rho_g)$, when $\omega=0.3$, $\eta=0.3$ and $M=1$: the black curve represents the mean value of 10 arbitrary samples.}
\label{Ch4_2dFilterFedEG}
\end{figure}

\begin{figure}[thpb]
\centering
\includegraphics[width=6cm,height=5cm]{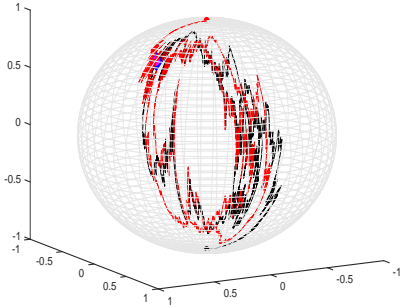}
\caption{\small Behavior of the trajectories $\rho_t$ and $\hat{\rho}_t$ with the feedback law $u(\hat{\rho})\equiv 1$, starting at $(\boldsymbol\rho_e,\boldsymbol\rho_g)$, when $\omega=0.3$, $\eta=0.3$ and $M=1$. The black curve corresponds to a sample trajectory $\rho_t$, the blue point represents its end point; the red curve corresponds to the  estimated trajectory $\hat{\rho}_t$, the magenta point represents its end point.}
\label{Ch4_2dFilter3DEG}
\end{figure}

In this section, we illustrate Proposition~\ref{prop:fidelity} through simulations of the system~\eqref{Ch5_StabFilter:2D SME}--\eqref{Ch5_StabFilter:2D SME filter}  in the case $u(\hat{\rho})\equiv 1$ and with parameters $\omega=0.3$, $\eta=0.3$ and $M=1.$ We set $\boldsymbol\rho_e$ as the initial state of the actual quantum filter and $\boldsymbol\rho_g$ as the initial state of the quantum filter esimate.
In Fig.~\ref{Ch4_2dFilterFedEG}, we simulate the fidelity $\mathcal{F}(\rho_t,\hat{\rho}_t)$.  Fig.~\ref{Ch4_2dFilter3DEG} represents the behavior of a sample trajectory $\rho_t$ and its corresponding estimation  $\hat{\rho}_t$ in Bloch sphere coordinates. 
\section{Asymptotics and feedback control of the coupled spin-$J$ systems}
\label{Ch5_Sec:ExpStabCoupledSys}
In this section, we discuss the asymptotic behavior of the actual quantum filter and its estimate for spin-$J$ systems with unknown initial states. The stochastic master equations are given by
\begin{align}
d\rho_t&\!=\!F_{\hat u_t}(\rho_t)dt\!+\!L(\rho_t)dt\!+\!G(\rho_t)\big(dY_t\!-\!2\sqrt{\eta M}\mathrm{Tr}(J_z \rho_t)dt\big),\nonumber\\
d\hat{\rho}_t&\!=\!F_{\hat u_t}(\hat{\rho}_t)dt\!+\!L(\hat{\rho}_t)dt\!+\!G(\hat{\rho}_t)\big(dY_t\!-\!2\sqrt{\eta M}\mathrm{Tr}(J_z \hat{\rho}_t)dt\big),\nonumber
\end{align}
where
\begin{itemize}
\item the actual quantum state of the spin system is denoted as $\rho$, and belongs to the space $
\mathcal{S}_N:=\{\rho\in\mathbb{C}^{N\times N}|\,\rho=\rho^*,\mathrm{Tr}(\rho)=1,\rho \geq 0\}$. The associated estimated state is denoted as $\hat{\rho}\in\mathcal{S}_N$,
\item $F_{u}(\rho):=-i[\omega J_z+u \,J_y,\rho]$, $L(\rho):=\frac{M}{2}(2J_z\rho J_z-J_z^2\rho-\rho J_z^2)$ and $G(\rho):=\sqrt{\eta M}\big(J_z\rho+\rho J_z-2\mathrm{Tr}(J_z\rho)\rho\big)$,
\item $Y_t$ denotes the observation process of the actual quantum spin system, which is a continuous semi-martingale whose quadratic variation is given by $[Y,Y]_t=t$. Its dynamics satisfies $dY_t=dW_t+2\sqrt{\eta M}\mathrm{Tr}(J_z\rho_t)dt$, where $W_t$ is a one-dimensional standard Wiener process,
\item $\hat u_t:=u(\hat{\rho}_t)$ denotes the feedback controller as a function of the estimated state $\hat{\rho}_t$,
\item 
$J_z$ is the (self-adjoint) angular momentum along the  $z$ axis, and it is defined by 
$
J_z e_n=(J-n)e_n,\quad n\in\{0,\dots,2J\}, 
$
where $J=\frac{N-1}{2}$ represents the fixed angular momentum and $\{e_0,\dots,e_{2J}\}$ corresponds to an orthonormal basis of $\mathbb C^N.$ With respect to this basis, the matrix form of $J_z$ is given by
\begin{equation}
J_z=
\begin{bmatrix}
J &&&&  \\
& J-1&&&\\
&&\ddots&&\\
&&&-J+1&\\
&&&&-J
\end{bmatrix},
\label{Ch3_Stabspin:Jz}
\end{equation}
We define the pure states $\boldsymbol \rho_n:=e_n e_n^*$ for $n\in\{0,\dots,2J\}$ corresponding to the eigenvectors of~$J_z$. 
\item 
$J_y$ is the (self-adjoint) angular momentum along the $y$ axis, and it is defined by 
\begin{equation}
J_ye_n=-ic_{n}e_{n-1}+ic_{n+1}e_{n+1},\quad n\in\{0,\dots,2J\},
\label{Ch3_Stabspin:Jy}
\end{equation}
where $c_m=\frac12\sqrt{(2J+1-m)m}$. The matrix form of $J_y$ is given by
\begin{equation*}
J_y=
\begin{bmatrix}
0&-ic_1 &&&\\
ic_1&0&-ic_2&&\\
&\ddots&\ddots&\ddots&\\
&&ic_{2J-1}&0&-ic_{2J}\\
&&&ic_{2J}&0
\end{bmatrix},
\end{equation*}
\item $\eta\in(0,1]$ measures the efficiency of the detectors, $M>0$ is the strength of the interaction between the system and the probe, and $\omega \geq 0$ is a parameter characterizing the free Hamiltonian. 
\end{itemize}
By replacing  $dY_t=dW_t+2\sqrt{\eta M}\mathrm{Tr}(J_z\rho_t)dt$ in the above equation, we obtain the following matrix-valued stochastic differential equations describing the time evolution of the pair $(\rho_t,\hat{\rho}_t)\in\mathcal{S}_N\times\mathcal{S}_N$,
\begin{align}
d\rho_t&\!=\!F_{\hat u_t}(\rho_t)dt+L(\rho_t)dt+G(\rho_t)dW_t,\label{Ch5_StabFilter:ND SME}\\
d\hat{\rho}_t&\!=\!F_{\hat u_t}(\hat{\rho}_t)dt+L(\hat{\rho}_t)dt+2\sqrt{\eta M}G(\hat{\rho}_t)\mathrm{Tr}\big(J_z(\rho_t-\hat{\rho}_t)\big)dt\nonumber\\
&\ +G(\hat{\rho}_t)dW_t.\label{Ch5_StabFilter:ND SME filter}
\end{align} 
If $u\in\mathcal{C}^1(\mathcal{S}_N,\mathbb{R})$, the existence and uniqueness of solutions of~\eqref{Ch5_StabFilter:ND SME}--\eqref{Ch5_StabFilter:ND SME filter} can be shown by similar arguments as in~\cite[Proposition 3.5]{mirrahimi2007stabilizing}.

Note that, if we turn off the feedback controller, there are $N^2$ equilibria $(\boldsymbol \rho_{n},\boldsymbol \rho_{m})$ with $n,m\in\{0,\dots,2J\}$ for the coupled system. However, since the system~\eqref{Ch5_StabFilter:ND SME} satisfies the non-demolition condition~\cite[Definition 2]{benoist2014large} and the measurement operator $J_z$ satisfies the non-degeneracy condition~\cite[Assumption (\textbf{ND})]{benoist2014large}, based on~\cite[Proposition 3]{benoist2014large}, we may state the following result.
\begin{theorem}[\cite{benoist2014large}]
 If $u\equiv 0$ and $\eta=1,$ $(\rho_t,\hat{\rho}_t)$ converges exponentially  towards the  set
$
\{(\boldsymbol \rho_{0},\boldsymbol \rho_{0}),\dots,(\boldsymbol \rho_{2J},\boldsymbol \rho_{2J})\}.
$
\label{thm:qsr}
\end{theorem}
In the following lemma, we will show that the fidelity is a sub-martingale. 
\begin{lemma}
Consider the coupled system~\eqref{Ch5_StabFilter:ND SME}--\eqref{Ch5_StabFilter:ND SME filter} with $\eta\in(0,1]$ and $u\in\mathcal{C}^1(\mathcal{S}_N,\mathbb{R}).$ Then the fidelity $\mathcal F(\rho,\hat\rho)$ is a sub-martingale.
\label{lem:sm}
\end{lemma}
Indeed, note that the value $\mathscr{L} \mathcal F(\rho,\hat{\rho})$ only depends on the instantaneous value of the Hamiltonian  $\omega J_z+u(\hat \rho) J_y$. Then by~\cite[Theorem~5]{amini2014stability}, one has $\mathscr{L} \mathcal F(\rho,\hat{\rho})\geq 0$ for the constant Hamiltonian $H=\omega J_z+u(\hat \rho) J_y$ which proves the result.
\begin{remark}
The infinitesimal generator $\mathscr{L} \mathcal F(\rho,\hat{\rho})$ can be decomposed in two terms, only one of which depends on the Hamiltonian. This term can be shown to be always zero. This is consistent with the above lemma.
\end{remark}
\begin{remark}
In Proposition~\ref{prop:fidelity}, by explicit computations, we have shown that for spin-$\frac 12$ systems $\mathscr{L} \mathcal F(\rho,\hat{\rho})\geq 0,$ which is consistent with Lemma~\ref{lem:sm}. This implied that $\mathscr{L} \mathcal F(\rho,\hat{\rho})$ should converge to zero almost surely. Moreover, we were able to analyze the zeros of $\mathscr{L} \mathcal F(\rho,\hat{\rho})$ and deduce some information about  the asymptotic behavior of $\mathcal F(\rho,\hat{\rho}).$ Similar results for spin-$J$ systems cannot be obtained easily because of the complicated form of  $\mathscr{L} \mathcal F(\rho,\hat{\rho}).$ 
\end{remark}
The above results give some intuitions about the asymptotic behavior of $\rho$ and $\hat\rho$ at least for the case $\eta=1$ or for spin-$\frac 12$ systems. 
 In these cases, without loss of generality, one may assume that the fidelity $\mathcal{F}(\rho_0,\hat\rho_0)$ is close to one. If $\eta=1,$ this can be obtained for example by turning off the feedback controller for a large enough time in view of Theorem~\ref{thm:qsr}. For spin-$\frac12$ systems, one can exploit Proposition~\ref{prop:fidelity} and apply a feedback controller satisfying $u(\boldsymbol\rho_{g})u(\boldsymbol\rho_{e})\neq 0$ for a large enough time. Now by Lemma~\ref{lem:sm}, we have $\mathbb E(\mathcal F(\rho_t,\hat\rho_t))\geq \mathcal F(\rho_0,\hat\rho_0)$ which leads to the conclusion that $\mathbb E(\mathcal F(\rho_t,\hat\rho_t))$ is close to one, independently of the chosen feedback law. 

In the following, we discuss the possibility of designing a feedback controller $u(\hat{\rho})$ which stabilizes exponentially almost surely the coupled system~\eqref{Ch5_StabFilter:ND SME}--\eqref{Ch5_StabFilter:ND SME filter} towards a given target state $(\boldsymbol \rho_{\bar{n}},\boldsymbol \rho_{\bar{n}})$ with $\bar{n}\in\{0,\dots,2J\}$. 
We note that the two subsystems~\eqref{Ch5_StabFilter:ND SME}--\eqref{Ch5_StabFilter:ND SME filter} share the same feedback controller $u(\hat{\rho})$, only depending on the estimated state $\hat{\rho}$. Hence, if we suppose that the feedback $u$ satisfies the assumption
$u(\boldsymbol\rho_{\bar{n}})=0$ and $u(\boldsymbol\rho_{k})\neq 0$ for all $k\neq \bar{n}$, then the coupled system~\eqref{Ch5_StabFilter:ND SME}--\eqref{Ch5_StabFilter:ND SME filter} possesses exactly $N$ equilibria, given by $(\boldsymbol\rho_{k},\boldsymbol\rho_{\bar{n}})$ for $k\in\{0,\dots,2J\}$.

In the aim of feedback exponential stabilization of the coupled spin-$J$ systems~\eqref{Ch5_StabFilter:ND SME}--\eqref{Ch5_StabFilter:ND SME filter} towards $(\boldsymbol \rho_{\bar{n}},\boldsymbol \rho_{\bar{n}})$, we propose a conjecture inspired by~\cite[Theorem 6.4 and Theorem 6.5]{liang2019exponential}. These results were developed for the case  $\rho_0=\hat\rho_0.$ 
\begin{conjecture}
  Consider the coupled system~\eqref{Ch5_StabFilter:ND SME}--\eqref{Ch5_StabFilter:ND SME filter} with $(\rho_0,\hat{\rho}_0)\in \mathcal S_N\times \mathcal S_N\setminus (\boldsymbol\rho_{\bar n},\boldsymbol\rho_{\bar n})$ and assume $\eta\in(0,1).$ Then, the feedback controller 
\begin{equation}
u_{\bar{n}}(\hat{\rho})  = \alpha \big(1-\mathrm{Tr}(\hat{\rho} \boldsymbol \rho_{\bar{n}})\big)^{\beta}, \qquad \alpha> 0,\quad \beta \geq 1,
\label{u_t 0 2J}
\end{equation}
exponentially stabilizes $(\rho_t,\hat{\rho}_t)$  to $(\boldsymbol \rho_{\bar{n}},\boldsymbol \rho_{\bar{n}})$ almost surely for the special case $\bar{n}\in\{0,2J\}$ with sample Lyapunov exponent less or equal than $-\eta M$. Moreover, the feedback controller
\begin{equation}
u_{\bar{n}}(\hat{\rho}) = \alpha \big(J-\bar{n}-\mathrm{Tr}(J_z\hat{\rho})\big)^{\beta}, \qquad \alpha> 0,\quad \beta \geq 1,
\label{u_t N}
\end{equation}
 exponentially stabilizes $(\rho_t,\hat{\rho}_t)$ to $(\boldsymbol \rho_{\bar{n}},\boldsymbol \rho_{\bar{n}})$ almost surely for the general case $\bar{n}\in\{0,\dots,2J\}$ with sample Lyapunov exponent less or equal than $-\eta M/2$ for $\bar{n}\in\{1,\dots,2J-1\}$ and $-\eta M$ for  $\bar{n}\in\{0,2J\}$. 
 \label{conj:fs}
\end{conjecture}
\subsection*{Simulations}
In this section, we illustrate the above conjecture through simulations for a three-level quantum spin system. 

\begin{figure}[thpb]
\centering
\includegraphics[width=9.7cm]{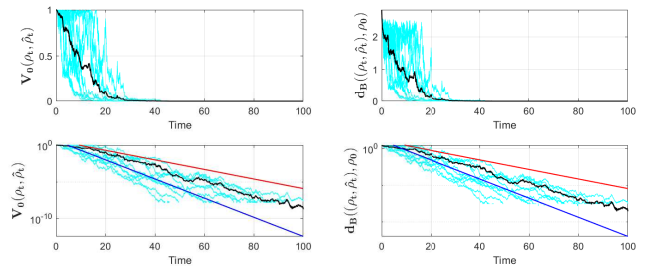}
\caption{\small Exponential stabilization of a three-level quantum spin system towards $(\boldsymbol\rho_0,\boldsymbol\rho_0)$ with the feedback law~\eqref{u_t 0 2J} starting at $(\rho_0,\hat{\rho}_0)=(\boldsymbol \rho_2,\boldsymbol \rho_1)$ with $\omega=0.3$, $\eta=0.3$, $M=1$, $\alpha = 5$ and $\beta = 2$. The black curve represents the mean value of 10 arbitrary sample trajectories, the red and blue curves represent the exponential references with exponents $-\eta M/2$ and $-\eta M$ respectively. The figures at the bottom are the semi-log versions of the ones at the top.}
\label{Ch5_3DQuantumFilter00}
\end{figure}

\begin{figure}[thpb]
\centering
\includegraphics[width=9.7cm]{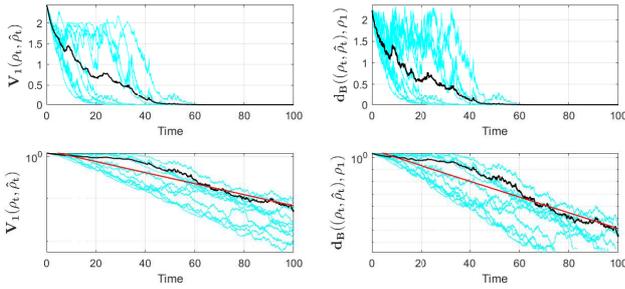}
\caption{\small Exponential stabilization of a three-level quantum spin system towards  $(\boldsymbol \rho_1,\boldsymbol \rho_1)$ with the feedback law~\eqref{u_t N} starting at $(\rho_0,\hat{\rho}_0)=\big(\mathrm{diag}(0.2,0.2,0.6),\mathrm{diag}(0.8,0.1,0.1)\big)$ with $\omega=0.3$, $\eta=0.3$, $M=1$, $\alpha = 2$ and $\beta = 2$. The black curve represents the mean value of 10 arbitrary sample trajectories, the red curve represents the exponential reference with exponent $-\eta M/2$. The figures at the bottom are the semi-log versions of the ones at the top.}
\label{Ch5_3DQuantumFilter11}
\end{figure}

In Fig~\ref{Ch5_3DQuantumFilter00}, we consider the target state $(\boldsymbol\rho_0,\boldsymbol\rho_0)$, and a candidate Lyapunov function $V_0(\rho,\hat{\rho})=\sqrt{1-\mathrm{Tr}(\rho\rho_{0})\mathrm{Tr}(\hat{\rho}\hat{\rho}_{0})}.$ We define the distance for the coupled system via the Bures distance  $d_B(\cdot,\cdot)$ in $S_3$~(\cite{bengtsson2017geometry}) as follows $\mathbf{d_B}((\rho,\hat{\rho}),\boldsymbol\rho_0):=d_B(\rho,\boldsymbol\rho_0)+d_B(\hat{\rho},\boldsymbol\rho_0)$. Then we have $\frac{\sqrt{2}}{4}V_0(\rho,\hat{\rho})\leq\mathbf{d_B}((\rho,\hat{\rho}),\boldsymbol\rho_0 )\leq \sqrt{2}V_0(\rho,\hat{\rho})$. 
In Fig~\ref{Ch5_3DQuantumFilter11}, we choose  $(\boldsymbol \rho_1,\boldsymbol \rho_1)$ as the target state and  $V_1(\rho,\hat{\rho})=\sum_{k\neq 1}\left(\sqrt{\mathrm{Tr}(\rho\rho_{k})}+\sqrt{\mathrm{Tr}(\hat\rho\hat{\rho}_{k})}\right)$ as candidate Lyapunov function. We can show that $\frac{\sqrt{2}}{2}V_1(\rho,\hat{\rho})\leq\mathbf{d_B}((\rho,\hat{\rho}),\boldsymbol\rho_1 )\leq \sqrt{2}V_1(\rho,\hat{\rho})$. 
\section{Conclusion}
In this paper, for spin-$\frac12$ systems, we have shown that under appropriate assumption on the feedback controller, we can guarantee the same asymptotics for the actual quantum filter and its estimate. For the general case of spin-$J$ systems, we provide an heuristic approach regarding the asymptotic behavior of $\rho$ and $\hat\rho$ and feedback exponential stabilization of the coupled system~\eqref{Ch5_StabFilter:ND SME}--\eqref{Ch5_StabFilter:ND SME filter}. In the future work, we aim to provide a rigorous proof of the conjecture~\ref{conj:fs} and present systematic methods ensuring exponential stabilization of the coupled spin systems. Other research lines include generalization to more general forms of filter equations which are driven by Wiener and/or Poisson processes. Also, exponential  stabilization of entangled states such as GHZ states, with unknown initial states is included in our research lines (see~\cite{liang2019exponential_two-qubit} for exponential stabilization of two-qubit systems with perfect initialization).

\bibliography{ref}

\end{document}